\documentclass[aos,MSNbibl,dvips]{arximspdf}
\usepackage{graphicx}

\doi{10.1214/13-AOS1113} 
\volume{41}
\issue{3}
\pubyear{2013}
\firstpage{1329}
\lastpage{1349}

\makeatletter
\newcommand{\rrvert}{\vert}
\newcommand{\llvert}{\vert}

\newproclaim{assumption}{Assumption}[section]

\newtheorem{theorem}{Theorem}[section]

\newproclaim{remark}{Remark}[section]

\newtheorem{lemma}{Lemma}[section]

\newtheorem{proposition}{Proposition}[section]

\def\var{\operatorname{var}}

\makeatother

\begin{document}
\begin{frontmatter}

\title{Fixed-smoothing asymptotics for time series\thanksref{T1}}
\runtitle{Fixed-smoothing asymptotics for time series}

\thankstext{T1}{Supported in part by NSF Grants DMS-08-04937 and DMS-11-04545.}

\begin{aug}
\author[A]{\fnms{Xianyang} \snm{Zhang}\corref{}\ead[label=e1]{zhang104@illinois.edu}}
\and
\author[B]{\fnms{Xiaofeng} \snm{Shao}\ead[label=e2]{xshao@illinois.edu}}

\runauthor{X. Zhang and X. Shao}

\affiliation{University of Missouri-Columbia and
University of Illinois at~Urbana-Champaign}

\address[A]{Department of Statistics\\
University of Missouri-Columbia\\
Columbia, Missouri 65211\\
USA\\
\printead{e1}}
\address[B]{Department of Statistics\\
University of Illinois\\
\quad at Urbana-Champaign\\
Champaign, Illinois 61820\\
USA\\
\printead{e2}} 
\end{aug}

\received{\smonth{12} \syear{2012}}
\revised{\smonth{3} \syear{2013}}

%
\begin{abstract}
In this paper, we derive higher order Edgeworth expansions for the
finite sample distributions of the subsampling-based $t$-statistic and
the Wald statistic in the Gaussian location model under the so-called
fixed-smoothing paradigm. In particular, we show that the error of
asymptotic approximation is at the order of the reciprocal of the
sample size and obtain explicit forms for the leading error terms in
the expansions. The results are used to justify the second-order
correctness of a new bootstrap method, the Gaussian dependent
bootstrap, in the context of Gaussian location model.
\end{abstract}

%

\end{frontmatter}

\section{Introduction}
\label{secintro}

Many economic and financial applications involve time series data
with autocorrelation and heteroskedasticity properties. Often the
unknown dependence structure is not the chief object of interest but
the inference on the parameter of interest involves the estimation
of unknown dependence. In stationary time series models estimated by
generalized method of moments (GMM), robust inference is typically
accomplished by consistently estimating the asymptotic covariance
matrix, which is proportional to the long run variance (LRV) matrix
of the estimating equations or moment conditions defining the
estimator, using a kernel smoothing method. In the econometrics and
statistics literature, the bandwidth parameter/truncation lag
involved in the kernel smoothing method is assumed to grow slowly
with sample size in order to achieve consistency. The inference is
conducted by plugging in a covariance matrix estimator that is
consistent under heteroskedasticity and autocorrelation. This
approach dates back to Newey and West~\cite{r24} and Andrews
\cite{r2}. Recently, Kiefer and Vogelsang~\cite{r15} (KV, hereafter)
developed an alternative first-order asymptotic theory for the HAC
(heteroskedasticity and autocorrelation consistent) based robust
inference, where the proportion of the bandwidth involved in the HAC
estimator to the sample size $T$, denoted as $b$, is held fixed in
the asymptotics. Under the fixed-$b$ asymptotics, the HAC estimator
converges to a nondegenerate yet nonstandard limiting distribution.
The tests based on the fixed-$b$ asymptotic approximation were shown
to enjoy better finite sample properties than the tests based on the
small-$b$ asymptotic theory under which the HAC estimator is
consistent, and the limiting distribution of the studentized
statistic admits a standard form, such as standard normal or
$\chi^2$ distribution. Using the higher order Edgeworth expansions,
Jansson~\cite{r14}, Sun et al.~\cite{r29} and Sun~\cite{r26}
rigorously proved that the fixed-$b$ asymptotics provides a high-order
refinement over the traditional small-$b$ asymptotics in the
Gaussian location model. Sun et al.~\cite{r29} also provided an
interesting decision theoretical justification for the use of
fixed-$b$ rules in econometric testing. For non-Gaussian linear
processes, Gon\c{c}alves and Vogelsang~\cite{r6} obtained an upper
bound on the convergence rate of the error in the fixed-$b$
approximation and showed that it can be smaller than the error of
the normal approximation under suitable assumptions.

Since the seminal contribution by KV, there has been a growing body
of work in econometrics and statistics to extend and expand the
fixed-$b$ idea in the inference for time series data. For example,
Sun~\cite{r27} developed a procedure for hypothesis testing in time
series models by using the nonparametric series method. The basic
idea is to project the time series onto a space spanned by a set of
fourier basis functions (see Phillips~\cite{r25} and M\"{u}ller
\cite{r235} for early developments) and construct the covariance
matrix estimator based on the projection vectors with the number of
basis functions held fixed. Also see Sun~\cite{r28} for the use of a
similar idea in the inference of the trend regression models.
Ibragimov and M\"{u}ller~\cite{r12} proposed a subsampling based
$t$-statistic for robust inference where the unknown dependence
structure can be in the temporal, spatial or other forms.
In their paper, the number of non-overlapping blocks is
held fixed. The $t$-statistic-based approach was extended by Bester
et al.~\cite{r3} to the inference of spatial and panel data with
group structure. In the context of misspecification testing, Chen
and Qu~\cite{r5} proposed a modified $M$ test of Kuan and Lee
\cite{r16} which involves dividing the full sample into several
recursive subsamples and constructing a normalization matrix based
on them. In the statistical literature, Shao~\cite{r255} developed
the self-normalized approach to inference for time series data that
uses an inconsistent LRV estimator based on recursive subsample
estimates. The self-normalized method is an extension of Lobato
\cite{r215} from the sample autocovariances to more general
approximately linear statistics, and it coincides with KVs fixed-$b$
approach in the inference of the mean of a stationary time series by
using the Bartlett kernel and letting $b=1$. Although the above
inference procedures are proposed in different settings and for
different problems and data structures, they share a common feature
in the sense that the underlying smoothing parameters in the
asymptotic covariance matrix estimators such as the number of basis
functions, the number of cluster groups and the number of recursive
subsamples, play a similar role as the bandwidth in the HAC
estimator. Throughout the paper, we shall call these asymptotics,
where the smoothing parameter (or function of smoothing parameter)
is held fixed, the fixed-smoothing asymptotics. In contrast, when
the smoothing parameter grows with respect to sample size, we use
the term increasing-domain asymptotics. At some places the terms
fixed-$K$ (or fixed-$b$) and increasing-$K$ (or small-$b$)
asymptotics are used to follow the convention in the literature.

In this article, we derive higher order expansions of the finite
sample distributions of the subsampling-based $t$-statistic and the
Wald statistic with HAC covariance estimator when the underlying
smoothing parameters are held fixed, under the framework of the
Gaussian location model. Specifically, we show that the error in the
rejection probability (ERP, hereafter) is of order $O(1/T)$ under
the fixed-smoothing asymptotics. Under the assumption that the
eigenfunctions of the kernel in the HAC estimator have zero mean and
other mild assumptions, we derive the leading error term of order
$O(1/T)$ under the fixed-smoothing framework. These results are
similar to those obtained under the fixed-$b$ asymptotics (see Sun
et al.~\cite{r29}), but are stronger in the sense that we are able
to derive the exact form of the leading error term with order
$O(1/T)$. The explicit form of the leading error term in the
approximation provides a clear theoretical explanation for the
empirical findings in the literature regarding the direction and
magnitude of size distortion for time series with various degrees of
dependence. To the best of our knowledge, this is the first time
that the leading error terms are made explicit through the higher
order Edgeworth expansion under the fixed-smoothing asymptotics. It
is also worth noting that our nonstandard argument differs from that
in Jansson~\cite{r14} and Sun et al.~\cite{r29}, and it may be of
independent theoretical interest and be useful for future follow-up
work.

Second, we propose a novel bootstrap method for time series, the
Gaussian dependent bootstrap, which is able to mimic the second-order
properties of the original time series and produces a Gaussian
bootstrap sample. For the Gaussian location model, we show that the
inference based on the Gaussian dependent bootstrap is more accurate
than the first-order approximation under the fixed-smoothing
asymptotics. This seems to be the first time a bootstrap method is
shown to be second-order correct under the fixed-smoothing
asymptotics; see Gon\c{c}alves and Vogelsang~\cite{r6} for a recent
attempt for the moving block bootstrap in the non-Gaussian setting.

We now introduce some notation. For a vector
$x=(x_1,x_2,\ldots,x_{q_0})\in\mathbb{R}^{q_0}$, we let
$\|x\|=(\sum^{q_0}_{i=1}x_i^2)^{1/2}$ be the Euclidean norm. For a
matrix $A=(a_{ij})_{i,j=1}^{q_0}\in\mathbb{R}^{q_0\times q_0}$,
denote by $\|A\|_2=\sup_{\|x\|=1}\|Ax\|$ the spectral norm and
$\|A\|_{\infty}=\max_{1\leq i,j\leq q_0}|a_{ij}|$ the max
norm.\vspace*{1pt}
Denote by $\lfloor a\rfloor$ the integer part of a real number $a$.
Let $L^2[0,1]$ be the space of square integrable functions on
$[0,1]$. Denote by $D[0,1]$ the space of functions on $[0,1]$ which
are right continuous and have left limits, endowed with the
Skorokhod topology; see Billingsley~\cite{r4}. Denote by
``$\Rightarrow$'' weak convergence in the $\mathbb{R}^{q_0}$-valued
function space $D^{q_0}[0, 1]$, where $q_0\in\mathbb{N}$. Denote by
``$\rightarrow^d$'' and ``$\rightarrow^p$'' convergence in
distribution and convergence in probability, respectively. The
notation $N(\mu,\Sigma)$ is used to denote the multivariate normal
distribution with mean $\mu$ and covariance $\Sigma$. Let
$\chi^2_{k}$ be a random variable following $\chi^2$ distribution
with $k$ degrees of freedom and $G_k$ be the corresponding
distribution function.

The layout of the paper is as follows. Section~\ref{sechighorder}
contains the higher order expansions of the finite sample distributions
of the subsampling $t$-statistic and the Wald statistic with HAC
estimator. We introduce the Gaussian dependent bootstrap and the
results about its second-order accuracy in Section~\ref{secbootstrap}.
Section~\ref{secconclusion} concludes. Technical details and
simulation results are gathered in the supplementary material~\cite{r33}.

\section{Higher order expansions}
\label{sechighorder} This paper is partially motivated by recent
studies on the ERP for the Gaussian location model by Jansson
\cite{r14} and Sun et al.~\cite{r29}, who showed that the ERP is of
order $O(1/T)$ under the fixed-$b$ asymptotics, which is smaller
than the ERP under the small-$b$ asymptotics. A natural question is
to what extent the ERP result can be extended to the recently
proposed fixed-smoothing based inference methods under the
fixed-smoothing asymptotics. Following Jansson~\cite{r14} and Sun et
al.~\cite{r29}, we focus on the inference of the mean of a
univariate stationary Gaussian time series or equivalently, a
Gaussian location model. We conjecture that the higher order terms
in the asymptotic expansion under the Gaussian assumption will also
show up in the general expansion without the Gaussian assumption.

\subsection{Expansion for the finite sample distribution of
subsampling-based $t$-statistic} \label{subsett-statistic} We first
investigate the Edgeworth expansion of the finite sample distribution
of subsampling-based $t$-statistic (Ibragimov and M\"{u}ller
\cite{r12}). Here we treat the subsampling-based $t$-statistic and
other cases separately, because the $t$-statistic corresponds to a
different choice of normalization factor (compare with the Wald
statistic in Section~\ref{subsetwald-statistic}). Given the
observations $\{X_1,X_2,\ldots,X_T\}$ from a Gaussian stationary time
series, we divide the sample into $K$ approximately equal sized groups
of consecutive observations. The observation $X_i$ is in the $j$th
group if and only if $i\in\mathcal{M}_j=\{s\in\mathbb{Z}\dvtx
(j-1)T/K<s\leq jT/K\}, j=1,2,\ldots,K$. Define the sample mean of the
$k$th group as
\[
\hat{\mu}_k=\frac{1}{|\mathcal{M}_k|}\sum_{i\in\mathcal
{M}_k}X_i,\qquad
k=1,2,\ldots,K,
\]
where $|\cdot|$ denotes the cardinality of a finite set. Let
$\hat{\mu}=(\hat{\mu}_1,\hat{\mu}_2,\ldots,\hat{\mu}_K)'$,
$\bar{\mu}_n=\frac{1}{K}\sum^{K}_{i=1}\hat{\mu}_i$ and
$S^2_n=\frac{1}{K-1}\sum^{K}_{i=1}(\hat{\mu}_i-\bar{\mu}_n)^2$. Then
the subsampling-based $t$-statistic for testing the null hypothesis
$H_{0}\dvtx  \mu=\mu_0$ versus the alternative $H_{a}\dvtx  \mu\neq\mu_0$
is given by
%
\begin{equation}
\label{sub-t} T_{K}=\frac{\sqrt{K}(\bar{\mu}_n-\mu_0)}{S_n}
=\frac{\sqrt {K}(\bar{\mu}_n-\mu_0)}{ \{(\sum^{K}_{i=1}(\hat
{\mu}_i-\bar{\mu}_n)^2)/({K-1})\}^{1/2}}.
\end{equation}
Our goal here is to develop an Edgeworth expansion of $P(|T_K|\leq
x)$ when $K$ is fixed and sample size $T\rightarrow\infty$.
It is not hard to see that the distribution of $T_K$ is symmetric, so
it is sufficient to consider $P(|T_K|\leq x)$ since
$P(T_K\leq x)=\frac{1+P(|T_K|\leq x)}{2}$ for any $x\geq0$. Denote
by $t_k$ a random variable following $t$ distribution with $k$
degrees of freedom. The following theorem gives the higher order
expansion under the Gaussian assumption.
%
\begin{theorem}\label{t-stat}
Assume that $\{X_i\}$ is a stationary Gaussian time series
satisfying that $\sum^{+\infty}_{h=-\infty}\gamma_X(h)>0$ and
$\sum^{+\infty}_{h=-\infty}h^2|\gamma_X(h)|<\infty$. Further suppose
that $|\mathcal{M}_1|=|\mathcal{M}_2|=\cdots=|\mathcal{M}_K|$ and
$K$ is fixed. Then under $H_{0}$, we have
%
\begin{equation}
\sup_{x\in[0,+\infty)}\bigl|P\bigl(|T_K|\leq x\bigr)-\Psi(x;K)\bigr|=O
\bigl(1/T^2\bigr),
\end{equation}
where $\Psi(x;K)=P(|t_{K-1}|\leq
x)-\frac{B}{2\sigma^{2}T}\Upsilon(x;K)$ with
\begin{eqnarray*}
\Upsilon(x;K)&=&-K^2P\bigl(|t_{K-1}|\leq x\bigr)+(K+1)E \biggl[
\chi^2_{K-1}G_{1} \biggl(\frac{\chi
^2_{K-1}x^2}{K-1} \biggr)
\biggr]
\\
&&{}-E \biggl[\chi^2_1G_{K-1} \biggl(
\frac{(K-1)\chi^2_1}{x^2} \biggr) \biggr]+1
\end{eqnarray*}
and $B=\sum^{+\infty}_{h=-\infty}|h|\gamma_X(h)$.
\end{theorem}

\begin{figure}

\includegraphics{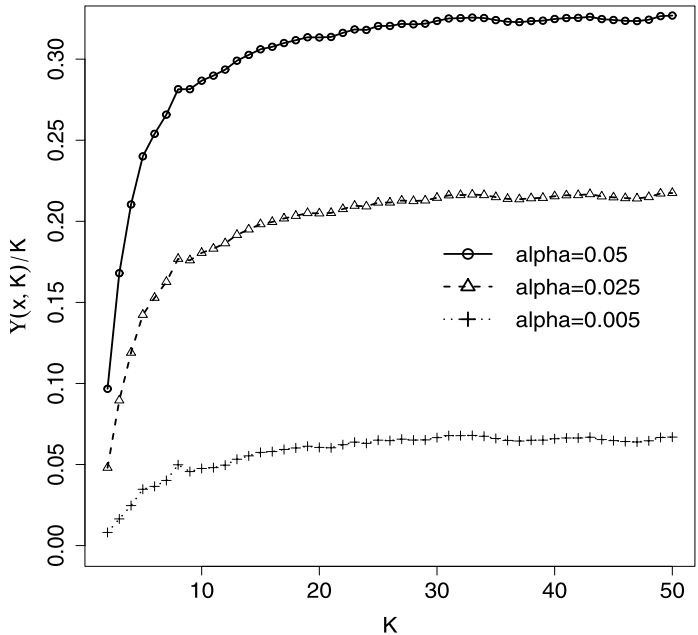}

\caption{Simulated values of $\Upsilon(t_{K-1}(1-\alpha);K)/K$ based
on 500,000 replications.}
\label{figtstat}
\end{figure}

We present the proof of Theorem~\ref{t-stat} in Section
\ref{proof-t-stat}, which requires some nonstandard arguments. From the
above expression, we see that the leading error term is of order
$O(1/T)$, and the magnitude and direction of the error depend upon
$B/\sigma^2$, which is related to the second-order properties of time
series, and $\Upsilon(x;K)$, which is independent of the dependence
structure of $\{X_i\}$ and can be approximated numerically for given
$x$ and $K$. Figure~\ref{figtstat} plots the approximated values of
$\Upsilon(t_{K-1}(1-\alpha);K)/K$ for different $K$ and $\alpha$, where
$t_{K-1}(1-\alpha)$ denotes the $100(1-\alpha)\%$ quantile of the $t$
distribution with $K-1$ degrees of freedom. It can be seen from
Figure~\ref{figtstat} that $\Upsilon(t_{K-1}(1-\alpha);K)/K$ increases
rapidly for $K<10$, and it becomes stable for relatively large $K$. For
each $K\geq2$, $\Upsilon(t_{K-1}(1-\alpha);K)/K$ is an increasing
function of $\alpha$. In the simulation work of Ibragimov and M\"uller
\cite{r12} (see Figure 2 therein), they found that the size of the
subsampling-based $t$-test is relatively robust to the correlations if
$K$ is small (say $K=4$ in their simulation). This finding is in fact
supported by our theory. For $K\le4$, the magnitude of $\Upsilon(x;K)$
is rather small, so the leading error term is small across a range of
correlations. As $K$ increases, the first-order approximation
deteriorates, which is reflected in the increasing magnitude of
$\Upsilon(t_{K-1}(1-\alpha);K)$ with respect to $K$.

Notice that $\Upsilon(t_{K-1}(1-\alpha);K)$ is always positive and
$\sigma^2>0$ by assumption, so the sign of the leading\vspace*{2pt} error term,
that is, $-\frac{B}{2\sigma^2T}\Upsilon(x;K)$, is determined by $B$.
When $B>0$ [e.g., AR(1) process with positive coefficient], the
first-order based inference tends to be oversized, and conversely it
tends to be undersized when \mbox{$B<0$} [e.g., MA(1) process with negative
coefficient]. Some simulations for AR(1) and MA(1) models in the
Gaussian location model support these theoretical findings. We
decide not to report these results to conserve space. Given the
sample size~$T$, the size distortion for the first-order based
inference may be severe if the ratio $B/\sigma^2$ is large. For
example, this is the case for AR(1) model, $X_t=\rho
X_{t-1}+\varepsilon_t$, as the correlation $\rho$ gets closer to 1.
As indicated by Figure~\ref{figtstat}, we show in the following
proposition that $\Upsilon(t_{K-1}(1-\alpha);K)/K$ converges as
$K\rightarrow+\infty$.
%
\begin{proposition}\label{leading-term}
As $K\rightarrow+\infty$, we have
$\Upsilon(x;K)/K=2x^2G_1'(x^2)+O(1/K)$, for any fixed
$x\in\mathbb{R}$.
\end{proposition}

Under the local alternative $H_{a}'\dvtx
\mu=\mu_0+(\delta\sigma)/\sqrt{T}$ with $\delta\neq0$, we can
derive a similar expansion for $T_K$ with $K$ fixed. Formally let
$Z$ be a random variable following the standard normal distribution
and $\mathcal{S}_{K-1}=\break\sqrt{\chi^2_{K-1}/(K-1)}$ with the
$\chi^2_{K-1}$ distribution being independent with $Z$. Then the
quantity $t_{K-1,\delta}=(Z+\delta)/\mathcal{S}_{K-1}$ follows a
noncentral $t$ distribution with noncentral parameter $\delta$.
Define $e_1(x)=E[\mathbf{I}\{|t_{K-1,\delta}|>x\}Z^2]$ and
$e_2(x)=E[\mathbf{I}\{|t_{K-1,\delta}|>x\}\chi^2_{K-1}]$. Then under
the local alternative, we have
\[
P\bigl(|T_K|\leq x\bigr)=P\bigl(|t_{K-1,\delta}|\leq x\bigr)-\frac{B}{2\sigma
^2T}
\Upsilon_{\delta}(x;K)+O\bigl(1/T^2\bigr),
\]
where
$\Upsilon_{\delta}(x;K)=K^2P(|t_{K-1,\delta}|>x)-e_1(x)-(K+1)e_2(x)$.
For fixed $\delta$, $P(|t_{K-1,\delta}|>t_{K-1}(1-\alpha))$ is a
monotonic increasing functions of $K$. An unreported numerical study
shows that $\Upsilon_{\delta}(t_{K-1}(1-\alpha);K)$ is roughly
monotonic with respect to $K$ for $\delta\in(0,4]$, which suggests
that larger $K$ tends to deliver more power when $B>0$. Combined
with the previous discussion, we see that the choice of $K$ leads to
a trade-off between the size distortion and power loss.

\begin{remark}
Theorem~\ref{t-stat} gives the ERP and the exact form of the
leading error term under the fixed-$K$ asymptotics. The higher order
expansion derived here is based on an expansion of the density
function of $(\hat{\mu}_1,\ldots,\hat{\mu}_K)$ which is made possible
by the Gaussian assumption. Extension to the general GMM setting
without the Gaussian assumption may require a different strategy in
the proof. Expansion for a distribution function or equivalently
characteristic function has been used in the higher order expansion
of the finite sample distribution under the Gaussian assumption (see,
e.g., Velasco and Robinson~\cite{r31} and Sun et al.~\cite{r29}).
With $K$ fixed in the asymptotics, the leading term of the variance
of the LRV estimator is captured by the first order fixed-$K$
limiting distribution and the leading term of the bias of the LRV
estimator is reflected in the leading error term
$-\frac{B}{2\sigma^2T}\Upsilon(x,K)$. Specifically, let
$\Sigma_T=(\sigma_{ij})^{K}_{i,j=1}$ with
$\sigma_{ij}=q\operatorname{Cov}(\hat{\mu}_i,\hat{\mu}_j)$. Then the leading
error term captures the difference between $\Sigma_T$ and
$\sigma^2I_K$, and the effect of the off-diagonal elements
$\sigma_{ij}$ with $|i-j|>1$ is of order $O(1/T^2)$ and thus is not
reflected in the leading term.
\end{remark}

\begin{remark}
When the number of groups $K$ grows slowly with the sample size
$T$, the Edgeworth expansion for $T_K$ was developed for $P(T_K\le
x)$ in Lahiri~\cite{r19,r20} under the general non-Gaussian setup.
The expansion given here is different from the usual Edgeworth
expansion under the increasing-domain asymptotics in terms of the
form and the convergence rate. Using the same argument, we can show
that under the fixed-$K$ asymptotics, the leading error term in the
expansion of $P(T_K\le x)$ is of order $O(1/T)$ under the Gaussian
assumption. In the non-Gaussian case, we conjecture that the order
of the leading error term is $O(1/\sqrt{T})$, which is due to the
effect of the third and fourth-order cumulants.
\end{remark}

The higher order Edgeworth expansion results in Sun et al.
\cite{r29} suggest that the fixed-$b$ based approximation is a
refinement of the approximation provided by the limiting
distribution derived under the small-$b$ asymptotics. In a similar
spirit, it is natural to ask if the fixed-$K$ based approximation
refines the first-order approximation under the increasing-$K$
asymptotics. To address this question, we consider the expansion
under the increasing-domain asymptotics, where $K$ grows slowly with
the sample size $T$.
%
\begin{proposition}\label{t-stat-increase}
Under the same conditions in Theorem~\ref{t-stat} but with
$\lim_{T\rightarrow\infty}(1/K+K/T)=0$, we have
%
\begin{equation}
\label{increasing}\quad P\bigl(|T_K|\leq x\bigr)=G_1
\bigl(x^2\bigr)+\frac{1}{K-1}x^4G_1''
\bigl(x^2\bigr)-\frac{BK}{T\sigma
^2}x^2G_1'
\bigl(x^2\bigr)+O(1/T).
\end{equation}
\end{proposition}
%
\begin{remark}
Since
\[
P\bigl(|t_{K-1}|\leq x\bigr)=G_1\bigl(x^2\bigr)+
\frac{1}{K-1}x^4G_1''
\bigl(x^2\bigr)+O\bigl(1/K^2\bigr)
\]
(see, e.g., Sun~\cite{r27}), we know that the fixed-$K$ based
approximation captures the first two terms in (\ref{increasing}),
whereas the increasing-$K$-based approximation (i.e., $\chi_1^2$)
only captures the first term. In view of Proposition
\ref{leading-term}, it is not hard to see that
\[
\Psi(x;K)=G_1\bigl(x^2\bigr)+\frac{1}{K-1}x^4G_1''
\bigl(x^2\bigr)-\frac{BK}{T\sigma
^2}x^2G_1'
\bigl(x^2\bigr)+O\bigl(1/K^2\bigr)+O(1/T),
\]
which implies that the fixed-$K$-based expansion is able to capture
all the three terms in (\ref{increasing}) as the smoothing parameter
$K\rightarrow\infty$ with $T^{1/3}=o(K)$. Loosely speaking, this
suggests that the fixed-$K$-based expansion holds for a broad range
of $K$, and it gets close to the corresponding increasing-$K$-based
expansion when $K$ is large.
\end{remark}

\subsection{Fixed-$b$ expansion}\label{subsetwald-statistic}
Consider a semi-positive definite bivariate kernel
$\mathcal{G}(\cdot,\cdot)$ which satisfies the spectral
decomposition
%
\begin{equation}
\label{spectral} \mathcal{G}(r,t)=\sum^{+\infty}_{j=1}
\lambda_j\phi_j(r)\phi _j(t),\qquad 0\leq r,t\leq1,
\end{equation}
where $\{\phi_j\}$ are the eigenfunctions, and $\{\lambda_j\}$ are
the eigenvalues which are in a descending order, that is,
$\lambda_1\geq\lambda_2\geq\cdots\geq0$. Suppose we have the
observations $\{X_1,X_2,\ldots,X_T\}$ from a stationary Gaussian time
series with mean $\mu$ and autocovariance function
$\gamma_X(i-j)=E[(X_i-\mu)(X_j-\mu)]$. The LRV estimator based on
the kernel $\mathcal{G}(\cdot,\cdot)$ and bandwidth $S_T=bT$ with
$b\in(0,1]$ is given by
\[
\hat{D}_{T,b}=\frac{1}{T}\sum_{i=1}^{T}
\sum_{j=1}^{T} \mathcal {G} \biggl(
\frac{i}{bT},\frac{j}{bT} \biggr) (X_i-
\bar{X}_T) (X_j-\bar{X}_T),
\]
where $\bar{X}_T=\sum^{T}_{i=1}X_i/T$ is the sample mean. For the
convenience of presentation, we set $b=1$. See Remark~\ref{remb}
for the case $b\in(0,1)$. To illustrate the idea, we define the
projection vectors
$\xi_j=\frac{1}{\sqrt{T}}\sum^{T}_{i=1}\phi_j^0(i/T)X_i$ with
$\phi_j^0(t)=\phi_j(t)-\frac{1}{T}\sum^{T}_{i=1}\phi_j(i/T)$ for
$j=1,2,\ldots\,$. Here\vspace*{1pt} the dependence of $\xi_j$ on $T$ is suppressed
to simplify the notation. Following Sun~\cite{r27}, we limit our
attention to the case $\int^{1}_{0}\phi_j(t)\,dt=0$ (e.g., Fourier
basis and Haar wavelet basis). For any semi-positive definite kernel
$\bar{\mathcal{G}}(\cdot,\cdot)$, we can define the demeaned kernel,
\[
\tilde{\mathcal{G}}(r,t)=\bar{\mathcal{G}}(r,t)-\int^1_0
\bar {\mathcal{G}}(s,t)\,ds-\int^1_0\bar{
\mathcal{G}}(r,p)\,dp+\int^1_0\int
^1_0\bar{\mathcal{G}}(s,p)\,ds\,dp.
\]
Suppose $\tilde{\mathcal{G}}(\cdot,\cdot)$ admits the spectral
decomposition
$\tilde{\mathcal{G}}(r,t)=\sum^{+\infty}_{i=1}\tilde{\lambda
}_i\tilde{\phi}_i(r)\tilde{\phi}_i(t)$
with $\{\tilde{\phi}_i\}$ and $\{\tilde{\lambda}_i\}$ being the
eigenfunctions and eigenvalues, respectively. Notice that
\[
\int^1_0\int^1_0
\tilde{\mathcal{G}}(r,t)\,dr\,dt=\sum^{+\infty}_{i=1}
\tilde{\lambda }_i \biggl(\int^1_0
\tilde{\phi}_i(t)\,dt \biggr)^2=0,
\]
which implies $\int^1_0\tilde{\phi}_i(t)\,dt=0$ whenever
$\lambda_i>0$, that is, the eigenfunctions of the demeaned kernel
$\tilde{\mathcal{G}}(\cdot,\cdot)$ are all mean zero. Based on the
spectral decomposition (\ref{spectral}) of
$\mathcal{G}(\cdot,\cdot)$, the LRV estimator with $b=1$ can be
rewritten as
\[
\hat{D}_{T,1}=\frac{1}{T}\sum_{i=1}^{T}
\sum_{j=1}^{T} \mathcal {G} \biggl(
\frac{i}{T},\frac{j}{T} \biggr) (X_i-
\bar{X}_T) (X_j-\bar {X}_T)=\sum
_{i=1}^{+\infty}\lambda_i\xi_i^2.
\]
We focus on testing the null hypothesis $H_{0}\dvtx  \mu=\mu_0$
versus the alternative $H_{a}\dvtx  \mu\neq\mu_0$. Define a sequence of
random variables
\[
F_T(K)=\frac{\xi_0^2}{\sum^{K}_{j=1}\lambda_j\xi_j^2},\qquad
K=1,\ldots,\infty,
\]
with $\xi_0=\frac{1}{\sqrt{T}}\sum^{T}_{i=1}(X_i-\mu_0)$. The Wald
test\vspace*{-1pt} statistic with HAC estiamtor is given by
$F_T(\infty)=\xi_0^2/\hat{D}_{T,1}$. Let $\{v_i\}_{i=0}^{+\infty}$
be a sequence of independent and identically distributed (i.i.d.)
standard normal\vadjust{\goodbreak} random variables. Further define
$\mathcal{F}(K):=\mathcal{F}(v;K)=\frac{v_0^2}{\sum^{K}_{j=1}\lambda_jv_j^2}$
and
%
\begin{eqnarray}
\label{psi} \aleph_{T}(x; K)=\frac{1}{2\sigma^2}\sum
^{K}_{i=0} \bigl(\operatorname{var}(\xi
_i)-\sigma^2 \bigr)E\bigl[\bigl(v_i^2-1
\bigr)\mathbf{I}\bigl\{\mathcal{F}(v;K)\leq
x\bigr\}\bigr],\nonumber\\[-8pt]\\[-8pt]
&&\eqntext{K=1,\ldots,\infty,}
\end{eqnarray}
with $\sigma^2=\sum^{+\infty}_{h=-\infty}\gamma_X(h)$ being the LRV.
The following theorem establishes the asymptotic expansion of the
finite sample distribution of $F_{T}(K)$ with $1\leq K\leq\infty$.
%
\begin{theorem}\label{fixed-b}
Assume the kernel $\mathcal{G}(\cdot,\cdot)$ satisfies the following
conditions:

(1) The second derivatives of the eigenfunctions $\{\phi
_i^{(2)}(\cdot)\}_{i=1}^{+\infty}$ exist. Further assume that the
eigenfunctions are mean zero and satisfy that
\[
\sup_{1\leq i\leq
J}\sup_{t\in[0,1]}\bigl|\phi_i^{(j)}(t)\bigr|<C J^j
\]
for $j=0,1,2$,
$J\in\mathbb{N}$, and some constant $C$ which does not depend on $j$~and~J;

(2) The eigenvalues $\lambda_n=O(1/n^{a})$, for some $a>19$.

Under the assumption that $\{X_i\}$ is a stationary Gaussian time series
with $\sigma^2=\sum^{+\infty}_{h=-\infty}\gamma_X(h)>0$ and
$\sum^{+\infty}_{h=-\infty}h^2|\gamma_X(h)|<\infty$, and the null
hypothesis $H_{0}$, we have $\sup_{x\in[0,+\infty)}|\aleph_{T}(x;
K)|=O(1/T)$ and
%
\begin{equation}
\label{second-order} \sup_{x\in[0,+\infty)}\bigl|P\bigl(F_T(K)\leq
x\bigr)-P\bigl(\mathcal{F}(K)\leq x\bigr)-\aleph_{T}(x; K)\bigr|=o(1/T)
\end{equation}
for any $1\leq K\leq\infty$.
\end{theorem}

The proof of Theorem~\ref{fixed-b} is based on the arguments of the
proof of Theorem~\ref{t-stat} given in Section~\ref{proof-t-stat}
and the truncation argument. The technical details are provided in
Zhang and Shao~\cite{r33}. For $K<\infty$, Theorem~\ref{fixed-b}
shows that the $O(1/T)$ ERP rate can be extended to the Wald
statistic with series variance estimator (Sun~\cite{r27}). When
$K=\infty$, Theorem~\ref{fixed-b} gives the asymptotic expansion of
the Wald test statistic $F_T(\infty)$ which is of particular
interest. The leading error term $\aleph_T(x;\infty)$ reflects the
departure of $\{\xi_j\}^{+\infty}_{j=0}$ from the i.i.d. standard
normal random variables $\{v_j\}^{+\infty}_{j=0}$. Specifically, the
form of $\aleph_T(x;\infty)$ suggests that the leading error term
captures the difference between the LRV and the variances of
$\xi_i$'s which are not exactly the same across $i=0,1,2,\ldots\,$. By
the orthogonality assumption, the covariance between $\xi_i$ and
$\xi_j$ with $i\neq j$ is of smaller order and hence is not
reflected in the leading term. Assume
$\int^{1}_{0}\mathcal{G}(r,r)\,dr=\sum^{+\infty}_{j=1}\lambda_j=1$. As
seen from Theorem~\ref{fixed-b}, the bias of the LRV estimator
[i.e., $\sum_{i=1}^{\infty} \lambda_i (\var(\xi_i)-\sigma^2)$] is
reflected in the leading error term $\aleph_{T}(x;\infty)$, which is\vadjust{\goodbreak}
a weighted sum of the relative difference of $\var(\xi_i)$ and
$\sigma^2$. Note that the difference $\var(\xi_i)-\sigma^2$ relies
on the second-order properties of the time series and the
eigenfunctions of $\mathcal{G}(\cdot,\cdot)$, and the weight
$E[(v_i^2-1)\mathbf{I}\{\mathcal{F}(\infty)\leq x\}]$ which depends
on the eigenvalues of $\mathcal{G}(\cdot,\cdot)$ is of order
$O(\lambda_i)$, as seen from the arguments used in the proof of
Theorem~\ref{fixed-b}.

In the econometrics and statistics literature, the bivariate kernel
$\mathcal{G}(\cdot,\cdot)$ is usually defined through a
semi-positive definite univariate kernel $\mathcal{K}(\cdot)$, that is,
$\mathcal{G}(r,t)=\mathcal{K}(r-t)$. In what follows, we make
several remarks regarding this special case.

\begin{remark}
\label{remb} For $0<b\leq1$, we define
$\mathcal{G}_b(\cdot,\cdot)=\mathcal{G}(\cdot/b,\cdot/b)$. If
$\mathcal{G}(\cdot,\cdot)$ is semi-positive definite on $[0,1/b]^2$,
then $\mathcal{G}_b(\cdot,\cdot)$ satisfies the spectral decomposition
$\mathcal{G}_b(r,t)=\sum^{+\infty}_{j=1}\lambda_{j,b}\phi
_{j,b}(r)\phi_{j,b}(t)$ with $0\leq r,t\leq1$. The eigencompoents of
$\mathcal{G}_b(r,t)$ can be obtained by solving a homogenuous Fredholm
integral equation of the second kind, where the solutions can be
approximated numerically when analytical solutions are unavailable.
When $\mathcal{G}(r,t)=\mathcal{K}(r-t)$, it was shown in Knessl and
Keller~\cite{r155} that under suitable assumptions on
$\mathcal{K}(\cdot)$,
$\lambda_{j,b}=b\int^{+\infty}_{-\infty}\mathcal{K}(r)\,dr-(\pi
^2j^2b^3/2)\int^{+\infty}_{-\infty}r^2\mathcal{K}(r)\,dr+o(b^3)$ and
$\phi_{j,b}\approx\sqrt{2}\sin(\pi jx)$ for $x$ bounded away from 0 and
1 as $b\rightarrow0$, which implies that
$\lambda_{M,b}/\lambda_{1,b}\rightarrow1$ for any fixed
$M\in\mathbb{N}$ and $b\rightarrow0$. Our result can be extended to the
case where $b<1$ if the assumptions in Theorem~\ref{fixed-b} hold for
$\{\lambda_{j,b}\}$ and $\{\phi_{j,b}\}$. It is also worth noting that
our result is established under different assumptions as compared to
Theorem 6 in Sun et al.~\cite{r29}, where the bivariate kernel is
defined as $\mathcal{G}(r,t)=\mathcal{K}(r-t)$ and the technical
assumption $b<1/(16\int^{+\infty}_{-\infty}|\mathcal{K}(r)|\,dr)$ is
required, which rules out the case $b=1$ for most kernels. Here we
provide an alternative way of proving the $O(1/T)$ ERP when the
eigenfunctions are mean zero. Furthermore, we provide the exact form of
the leading error term which has not been obtained in the literature.
\end{remark}

\begin{remark}
The assumption on the eigenvalues is satisfied by the bivariate
kernel defined through the QS kernel and the Daniel kernel with
$0<b\leq1$, and the Tukey--Hanning kernel with $b=1$ because these
kernels are analytical on the corresponding regions, and their
eigenvalues decay exponentially fast; see Little and Reade
\cite{r21}. However, the assumption does not hold for the Bartlett
kernel because the decay rate of its eigenvalues is of order
$O(1/n^2)$. For the demeaned Tukey--Hanning kernel with $b=1$, we
have that the eigenfunctions $\phi_1(t)=\sqrt{2}\cos\pi t$ and
$\phi_2(t)=\frac{\sin\pi t-2/\pi}{\sqrt{1/2-4/\pi^2}}$ with
eigenvalues $\lambda_1=0.25$, $\lambda_2=0.0474$ and $\lambda_j=0$
for $j\geq3$. It is not hard to construct a kernel that satisfies
the conditions in Theorem~\ref{fixed-b}. For example, one can
consider the kernel
$\mathcal{K}(r-t)=\sum^{+\infty}_{j=1}\lambda_j\{\cos(2\pi
jr)\cos(2\pi jt)+\sin(2\pi jr)\sin(2\pi
jt)\}=\sum^{+\infty}_{j=1}\lambda_j\cos(2\pi j(r-t))$ with
$\sum^{+\infty}_{j=1}\lambda_j=1$ and
$\lambda_j=O(1/j^{19+\epsilon})$ for some $\epsilon>0$. Then the
asymptotic expansion (\ref{second-order}) holds for the Wald
statistic based on the difference kernel
$\mathcal{G}(r,t)=\mathcal{K}(r-t)$.
\end{remark}

Define the Parzen characteristic exponent
\[
q=\max \biggl\{q_0\dvtx q_0\in\mathbb{Z}^+, g_{q_0}=
\lim_{x\rightarrow
0}\frac{1-\mathcal{K}(x)}{|x|^{q_0}}<\infty \biggr\}.
\]
For the Bartlett kernel $q$ is 1; For the Parzen and QS kernels, $q$
is equal to~2. Let $c_1=\int^{+\infty}_{-\infty}\mathcal{K}(x)\,dx$
and $c_2=\int^{+\infty}_{-\infty}\mathcal{K}^2(x)\,dx$. Further define
$\mathcal{F}_b(\infty)$ and $\aleph_{T,b}(x;\infty)$ with $\phi_j$
and $\lambda_j$ being replaced with $\phi_{j,b}$ and $\lambda_{j,b}$
in the definition of $\mathcal{F}(\infty)$ and
$\aleph_{T}(x;\infty)$. We summarize the first and second-order
%
\begin{table}
\caption{Asymptotic comparison between the first and
second-order approximations based on fixed-$b$ and small-$b$
asymptotics}\label{comp}
\begin{tabular*}{\tablewidth}{@{\extracolsep{\fill}}l c c@{}}
\hline
\textbf{Asymptotics} & \textbf{First order} & \textbf{Second order} \\
\hline
Fixed-$b$ & $P (\mathcal{F}_b(\infty)\leq x )$ & $P
(\mathcal{F}_b(\infty)\leq x )+\aleph_{T,b}(x;\infty)$ \\[2pt]
Small-$b$ & $G_1(x)$ &
\multicolumn{1}{l@{}}{$G_1(x)+(c_2G_1''(x)x^2-c_1G_1'(x)x)b$}\\[2pt]
&& \multicolumn{1}{l@{}}{$\qquad{}-\frac
{g_q\sum^{+\infty}_{h=-\infty}|h|^q\gamma_X(h)}{\sigma
^2(bT)^q}G_1'(x)x$}\\
\hline
\end{tabular*}
\end{table}
approximations for the distribution of studentized sample mean in the
Gaussian location model based on both fixed-$b$ and small-$b$
asymptotics in Table~\ref{comp} above. The formulas for the
second-order approximation under the small-$b$ asymptotics is from
Velasco and Robinson~\cite{r31}.
%
\begin{remark}
A few remarks are in order regarding Table~\ref{comp}. First of
all, it is worth noting that $P (\mathcal{F}_b(\infty)\leq
x )=G_1(x)+(c_2G_1''(x)x^2-c_1G_1'(x)x)b+O(b^2)$ as
$b\rightarrow0$ in Sun et al.~\cite{r29}, which suggests that the
fixed-$b$ limiting distribution captures the first two terms in the
higher order asymptotic expansion under the small-$b$ asymptotics
and thus provides a better approximation than the $\chi^2_1$
approximation. Second, it is interesting to compare the second-order
asymptotic expansions under the fixed-$b$ asymptotics and
small-$b$ asymptotics. We show in Proportion~\ref{fix-small} that
the higher order expansion under fixed-$b$ asymptotics is consistent
with the corresponding higher order expansion under small-$b$
asymptotics as $b$ approaches zero.
\end{remark}

Because our fixed-$b$ expansion is established under the assumption
that the eigenfunctions have mean zero, we shall consider the Wald
statistic $F_T(\infty)$ based on the demeaned kernel
$\tilde{\mathcal{G}}_b(r,t)=\mathcal{K}_b(r-t)-\int^1_0\mathcal
{K}_b(s-t)\,ds-\int^1_0\mathcal{K}_b(r-p)\,dp+\int^1_0\int^1_0\mathcal
{K}_b(s-p)\,ds\,dp$
with $\mathcal{K}_b(\cdot)=\mathcal{K}(\cdot/b)$ and $b\in(0,1]$.
Let $\{\tilde{\phi}_{j,b}\}$ and $\{\tilde{\lambda}_{j,b}\}$ be the
corresponding eigenfunctions and eigenvalues of
$\tilde{\mathcal{G}}_b(\cdot,\cdot)$.

\begin{proposition}\label{fix-small}
Suppose $\mathcal{K}(\cdot)\dvtx  \mathbb{R}\rightarrow[0,1]$ is
symmetric, semi-positive definite, piecewise smooth with
$\mathcal{K}(0)=1$ and $\int^{+\infty}_{0}x\mathcal{K}(x)\,dx<\infty$.
The Parzen characteristic exponent of $\mathcal{K}$ is no less than
one. Further assume that
%
\begin{equation}
\label{sup} \sup_{k\in\mathbb{N}}\Biggl\llvert \sum
_{i=1}^k\tilde{\lambda }_{i,b}\bigl(
\operatorname{var}(\tilde{\xi}_{i,b})-\sigma^2\bigr)\Biggr
\rrvert =O \Biggl(\sum_{i=1}^{+\infty}\tilde{
\lambda}_{i,b}\bigl(\operatorname{var}(\tilde {\xi}_{i,b})-
\sigma^2\bigr) \Biggr)
\end{equation}
as\vspace*{1pt} $b+1/(bT)\rightarrow0$, where $\tilde{\xi}_{i,b}$ is defined by
replacing $\phi_j$ with $\tilde{\phi}_{j,b}$ in the definition of
$\xi_i$. Then under the assumption that $\sigma^2>0$ and
$\sum^{+\infty}_{h=-\infty}h^2|\gamma_X(h)|<\infty$, we have
\[
\aleph_{T,b}(x;\infty)= -\frac{g_q\sum^{+\infty}_{h=-\infty}|h|^q\gamma_X(h)}{\sigma
^2(bT)^q}G_1'(x)x
\bigl(1+o(1)\bigr)+O(1/T)
\]
for fixed $x\in\mathbb{R}$, as $b\rightarrow0$ and $bT\rightarrow
+\infty$.
\end{proposition}
In Proposition~\ref{fix-small}, condition (\ref{sup}) is not
primitive, and it requires that the~bias for the LRV estimators based
on the kernel
$\tilde{\mathcal{G}}_{k,b}(r,t)=\break\sum^{k}_{i=1}\tilde{\lambda
}_{j,b}\tilde{\phi}_{j,b}(r)\tilde{\phi}_{j,b}(t)$
is at the same or smaller order of the bias for the LRV estimator
based on $\tilde{\mathcal{G}}_{b}(r,t)$. This condition simplifies
our technical arguments and it can be verified through a
case-by-case study. As shown in Proposition~\ref{fix-small}, the
fixed-$b$ expansion is consistent with the small-$b$ expansion as
$b$ approaches zero, and it is expected to be more accurate in terms
of approximating the finite sample distribution when $b$ is
relatively large. Overall speaking, the above result suggests that
the fixed-$b$ expansion provides a good approximation to the finite
sample distribution which holds for a broad range of $b$.

\section{Gaussian dependent bootstrap}
\label{secbootstrap} Given the higher order expansions presented in
Section~\ref{sechighorder}, it seems natural to investigate if
bootstrapping can help to improve the first-order approximation. Though
the higher order corrected critical values can also be obtained by
direct estimation of the leading error term, it involves estimation of
the eigencomponents of the kernel function and a choice of truncation
number for the leading error term $\aleph_T(x;\infty)$ [see
(\ref{psi})] besides estimating the second-order properties of the time
series. Therefore it is rather inconvenient to implement this
analytical approach because numerical or analytical calculation of
eigencomponents can be quite involved, the truncation number and the
bandwidth parameter used in estimating second-order properties are both
user-chosen numbers, and it seems difficult to come up with good rules
about their (optimal) choice in the current context. By contrast, the
bootstrap procedure proposed below, which involves only one user-chosen
number, aims to estimate the leading error term in an automatic fashion
and the computational cost is moderate given current high computing
power.

To present the idea, we again limit our attention to the univariate
Gaussian location model. Consider a consistent estimate of the
covariance matrix of $\{X_i\}^{T}_{i=1}$ which takes the form
$\hat{\Xi}(\omega;l)\in\mathbb{R}^{T\times T}$ with the $(i,j)$th
element given by $\omega_l(i-j)\hat{\gamma}_{X}(|i-j|)$ for
$i,j=1,2,\ldots,T$, where $\omega$ is a kernel function with
$\omega_l(\cdot)=\omega(\cdot/l)$ and
$\hat{\gamma}_{X}(h)=\frac{1}{T}\sum^{T-h}_{i=1}(X_i-\bar
{X}_T)(X_{i+h}-\bar{X}_T)$ for $h=0,1,2,\ldots,T-1$. Estimating the
covariance matrix of a stationary time series has been investigated by
a few researchers. See Wu and Pourahmadi~\cite{r32} for the use of a
banded sample covariance matrix and McMurry and Politis~\cite{r23} for
a tapered version of the sample covariance matrix. In what follows, we
shall consider the Bartlett kernel, that is,
$\omega(x)=(1-|x|)\mathbf{I}\{|x|<1\}$,\vspace*{1pt} which guarantees
to yield a semi-positive definite estimates, that is,
$\hat{\Xi}(\omega;l)\geq0$.

We now introduce a simple bootstrap procedure which can be shown to
be second-order correct. Suppose $X_1^*,\ldots,X_T^*$ is the
bootstrap sample generated from $N(0,\hat{\Xi}(\omega,l))$. It is
easy to see that $X_i^*$'s are stationary and Gaussian conditional
on the data. This is why we name this bootstrap method ``Gaussian
dependent bootstrap.'' There is a large literature on bootstrap for
time series; see Lahiri~\cite{r18} for a review. However, most of
the existing bootstrap methods do not deliver a conditionally
normally distributed bootstrap sample. Since our higher order
results are obtained under the Gaussian assumption, we need to
generate Gaussian bootstrap samples in order for our expansion
results to be useful.

Denote by $T_K^*$ the bootstrapped subsampling $t$-statistic
obtained by replacing $(X_1{-}\mu_0,X_2{-}\mu_0,\ldots,X_T{-}\mu_0)$ with
$(X^*_1,X^*_2,\ldots,X_T^*)$. Define the~bootstrapped projection
vectors $\xi_0^*=\frac{1}{\sqrt{T}}\sum^{T}_{j=1}X_j^*$ and
$\xi_j^*=\frac{1}{\sqrt{T}}\sum^{T}_{i=1}\phi^0_j(i/T)X_i^*$ for
$j=1,\ldots\,$. Let $P^*$ be the bootstrap probability measure
conditional on the data. The following theorems state the second-order
accuracy of the Gaussian dependent bootstrap in the univariate
Gaussian location model.

\begin{theorem}\label{bootstrap-fixed-K}
For the Gaussian location model, under the same conditions in
Theorem~\ref{t-stat} and $1/l+l^3/T\rightarrow0$, we have
%
\begin{equation}
\sup_{x\in[0,+\infty)}\bigl\llvert P\bigl(|T_K|\leq x\bigr)-P^*
\bigl(|T_K^*|\leq x\bigr)\bigr\rrvert =o_p(1/T).
\end{equation}
\end{theorem}

\begin{theorem}\label{bootstrap-fixed-b}
For the Gaussian location model, under the assumptions in Theorem
\ref{fixed-b} and that $1/l+l^3/T\rightarrow0$, we have
%
\begin{equation}
\sup_{x\in[0,+\infty)}\bigl\llvert P\bigl(F_{T}(\infty)\leq x
\bigr)-P^*\bigl(F_{T}^*(\infty)\leq x\bigr)\bigr\rrvert
=o_p(1/T),
\end{equation}
where
$F_T^*(\infty)=\frac{(\xi_0^*)^2}{\sum^{+\infty}_{j=1}\lambda
_j(\xi_j^*)^2}$
with $\{\lambda_j\}_{j=1}^{+\infty}$ given in (\ref{spectral}). Note
that $F_T^*(\infty)=(\xi_0^*)^2/\hat{D}_{T,1}^*$, where
$\hat{D}_{T,1}^*=T^{-1}\sum_{i,j=1}^{T}\mathcal
{G}(i/T,j/T)(X_i^*-\bar{X}_T^*)(X_j^*-\bar{X}_T^*)$
and $\bar{X}_T^*$ is the bootstrap sample mean.
\end{theorem}

\begin{remark}
The higher order terms in the small-$b$ expansion and the
increasing-$K$ expansion (see Table~\ref{comp} and Proposition
\ref{t-stat-increase}) depend on the second-order properties only
through the quantities $\sum^{+\infty}_{h=-\infty}|h|^k\gamma_X(h)$
for $k=0,1,\ldots,q$. It suggests that the Gaussian dependent
bootstrap also preserves the second-order accuracy under the
increasing-domain asymptotics provided that
\[
\sum^{+\infty}_{h=-\infty}|h|^{q+1}\gamma_{X}(h)<\infty.
\]
A rigorous proof is omitted due to space limitation.
\end{remark}

The bootstrap-based autocorrelation robust testing procedures have
been well studied in both econometrics and statistics literature
under the increasing-domain asymptotics. In the statistical
literature, Lahiri~\cite{r17} showed that for the studentized
$M$-estimator, the ERP of the moving block bootstrap (MBB)-based
one-sided testing procedure is of order $o_p(T^{-1/2})$ which
provides an asymptotic refinement to the normal approximation. Under
the framework of the smooth function model, G\"{o}tze and K\"{u}nsch
\cite{r7} showed that the ERP for the MBB-based one-sided test is of
order $O_p(T^{-3/4+\epsilon})$ for any $\epsilon>0$ when the HAC
estimator is constructed using the truncated kernel. Note that in
the latter paper, the HAC estimator used in the studentized
bootstrap statistic needs to take a different form from the original
HAC estimator to achieve the higher order accuracy. Also see Lahiri
\cite{r19} for a recent contribution. In the econometric literature,
the Edgeworth analysis for the block bootstrap has been conducted by
Hall and Horowitz~\cite{r95}, Andrews~\cite{rr25} and Inoue and
Shintani~\cite{r13}, among others, in the GMM framework. Within the
increasing-domain asymptotic framework, it is still unknown whether
the bootstrap can achieve an ERP of $o_p(1/T)$ when a HAC covariance
matrix estimator is used for studentization; see H\"{a}rdle,
Horowitz and Kreiss~\cite{r10}. Note that Hall and Horowitz
\cite{r95} and Andrews~\cite{rr25} obtained the $o_p(1/T)$ results
for symmetrical tests but they assumed the uncorrelatedness of the
moment conditions after finite lags. Note that all the above results
were obtained under the non-Gaussian assumption.

Within the fixed-smoothing asymptotic framework, Jansson~\cite{r14}
established that the error of the fixed-$b$ approximation to the
distribution of two-sided test statistic is of order $O(\log(T)/T)$
for the Gaussian location model and the case $b=1$, which was
further refined by Sun et al.~\cite{r29} by dropping the $\log(T)$
term. In the non-Gaussian setting, Gon\c{c}alves and Vogelsang
\cite{r6} showed that the fixed-$b$ approximation to the
distribution of one-sided test statistic has an ERP of order
$o(T^{-1/2+\epsilon})$ for any $\epsilon>0$ when all moments exist.
The latter authors further showed that the MBB (with i.i.d. bootstrap
as a special case) is able to replicate the fixed-$b$ limiting
distribution and thus provides more accurate approximation than the
normal approximation. However, because the exact form of the leading
error term was not obtained in their studies, their results seem not
directly applicable to show the higher order accuracy of bootstrap
under the fixed-$b$ asymptotics. Using the asymptotic expansion
results developed in Section~\ref{sechighorder}, we show that the
Gaussian dependent bootstrap can achieve an ERP of order $o_p(1/T)$
under the Gaussian assumption. This appears to be the first result
that shows the higher order accuracy of bootstrap under the
fixed-smoothing asymptotics. Our result also provides a positive
answer to the open question mentioned in H\"{a}rdle, Horowitz and
Kreiss~\cite{r10} that whether the bootstrap can achieve an ERP of
$o_p(1/T)$ in the dependence case when a HAC covariance matrix
estimator is used for studentization. It is worth noting that our
result is established for the symmetrical distribution functions
under the fixed-smoothing asymptotics and the Gaussian assumption.
It seems that in general the ERP of order $o_p(1/T)$ cannot be
achieved under the increasing-domain asymptotics or for the
non-Gaussian case. In the supplementary material~\cite{r33}, we provide some
simulation results which demonstrate the effectiveness of the
proposed Gaussian dependent bootstrap in both Gaussian and
non-Gaussian settings. The MBB is expected to be second-order
accurate, as seen from its empirical performance, but a rigorous
theoretical justification seems very difficult. Finally, we mention
that it is an important problem to choose $l$. For a given
criterion, the optimal $l$ presumably depends on the second-order
property of the time series in a sophisticated fashion. Some of the
rules proposed for block-based bootstrap (see Lahiri~\cite{r18},
Chapter 7) may still work, but a serious investigation is beyond the
scope of this article.

\section{Conclusion}
\label{secconclusion} In this paper, we derive the Edgeworth expansions
of the subsampling-based $t$-statistic and the Wald statistic with HAC
estimator in the Gaussian location model. Our work differs from the
existing ones in two important aspects: (i) the expansion is derived
under the fixed-smoothing asymptotics and the ERP of order $O(1/T)$ is
shown for a broad class of fixed-smoothing inference procedures; (ii)
we obtain an explicit form for the leading error term, which is
unavailable in the literature. An in-depth analysis of the behavior of
the leading error term when the smoothing parameter grows with sample
size (i.e., $K\rightarrow\infty$ in the subsampling $t$-statistic or
$b\rightarrow0$ in the Wald statistic with the HAC estimator) shows the
consistency of our results with the expansion results under the
increasing-domain asymptotics. Building on these expansions, we further
propose a new bootstrap method, the Gaussian dependent bootstrap, which
provides a higher order correction than the first-order fixed-smoothing
approximation.

We mention a few directions that are worthy of future research.
First, it would be interesting to relax the Gaussian assumption in
all the expansions we obtained in the paper. For non-Gaussian time
series, Edgeworth expansions have been obtained by G\"{o}tze and
Kunsch~\cite{r7}, Lahiri~\cite{r19,r20}, among others, for
studentized statistics of a smooth function model under weak
dependence assumption, but their results were derived under the
increasing-smoothing asymptotics. For the location model and
studentized sample mean, the extension to the non-Guassian case may
require an expansion of the corresponding characteristic function,
which involves calculation of the high-order cumulants under the
fixed-smoothing asymptotics. The detailed calculation of the high-order
terms can be quite involved and challenging. We conjecture
that under the fixed-smoothing asymptotics, the leading error term
in the expansion of its distribution function involves the third and
fourth-order cumulants, which reflects the non-Gaussianness, and the
order of the leading error term is $O(T^{-1/2})$ instead of
$O(T^{-1})$. Second, we expect that our expansion results will be
useful in the optimal choice of the smoothing parameter, the kernel
and its corresponding eigenvalues and eigenfunctions, for a given
loss function. The optimal choice of the smoothing parameter has
been addressed in Sun et al.~\cite{r29} using the expansion derived
under the increasing-smoothing asymptotics. As the finite sample
distribution is better approximated by the corresponding
fixed-smoothing based approximations at either first or second order
than its increasing-smoothing counterparts, the fixed-smoothing
asymptotic theory proves to be more relevant in terms of explaining
the finite sample results; see Gon\c{c}alves and Vogelsang
\cite{r6}. Therefore, it might be worth reconsidering the choice of
the optimal smoothing parameter under the fixed-smoothing
asymptotics. Third, we restrict our attention to the Gaussian
location model when deriving the higher order expansions. It would
be interesting to extend the results to the general GMM setting. A
recent attempt by Sun~\cite{r26} for the HAC-based inference seems
to suggest this is feasible. Finally, under the fixed-smoothing
asymptotics, the second correctness of the moving block bootstrap
for studentized sample mean, although suggested by the simulation
results, is still an open but challenging topic for future research.

\section{\texorpdfstring{Proof of Theorem \protect\ref{t-stat}}
{Proof of Theorem 2.1}}\label{proof-t-stat}
Consider the $K+1$-dimensional multivariate normal density function
which takes the form
\[
f(y,\Sigma)=(2\pi)^{-({K+1})/{2}}|\Sigma|^{-{1}/{2}}\exp \bigl(-
\tfrac{1}{2}y'\Sigma^{-1}y \bigr).
\]
We assume the $(i,j)$th element and the $(j,i)$th element of
$\Sigma$ are functionally unrelated. The results can be extended to
the case where symmetric matrix elements are considered functionally
equal; see, for example, McCulloch~\cite{r22}. In the following, we use
$\otimes$ to denote the Kronecker product in matrix algebra and use
$\mathrm{vec}$ to denote the operator that transforms a matrix into a
column vector\vspace*{1pt} by stacking the columns of the matrix one underneath
the other. For a vector $y\in\mathbb{R}^{l\times1}$ whose elements
are differential functions of a vector $x\in\mathbb{R}^{k\times1}$,
we define $\frac{\partial y}{\partial x}$ to be a $k\times l$ matrix
with the $(i,j)$th element being $\frac{\partial y_j}{\partial
x_i}$. The notation $u \asymp v$ represents $u=O(v)$ and $v=O(u)$.
We first present the following lemmas whose proofs are given in the
online supplement~\cite{r33}.

\begin{lemma}\label{first-deri}
\[
\frac{\partial
f}{\partial\operatorname{vec}(\Sigma)}(y,\Sigma)=\frac{f(y,\Sigma)}{2}\bigl\{ \bigl(
\Sigma^{-1}y\bigr)\otimes\bigl(\Sigma^{-1}y\bigr)-
\operatorname{vec}\bigl(\Sigma^{-1}\bigr)\bigr\}.
\]
\end{lemma}

\begin{lemma}\label{second-deri}
\begin{eqnarray*}
&&\frac{\partial^2 f}{\partial
\operatorname{vec}(\Sigma)\operatorname{vec}(\Sigma)}(y,\Sigma)
\\
&&\qquad=\frac{1}{4}\bigl\{\bigl(\Sigma^{-1}y\bigr)\otimes\bigl(
\Sigma^{-1}y\bigr)- \operatorname{vec}\bigl(\Sigma^{-1}\bigr)
\bigr\} \bigl\{\bigl(\Sigma^{-1}y\bigr)\otimes\bigl(
\Sigma^{-1}y\bigr)- \operatorname{vec}\bigl(\Sigma^{-1}\bigr)
\bigr\}'\\
&&\qquad\quad{}\times f(y,\Sigma)
\\
&&\qquad\quad{}-\frac{1}{2}\bigl\{\bigl(\Sigma^{-1}yy'
\Sigma^{-1}\bigr)\otimes\Sigma ^{-1}+\Sigma^{-1}
\otimes\bigl(\Sigma^{-1}yy'\Sigma^{-1}\bigr)-
\Sigma ^{-1}\otimes\Sigma^{-1}\bigr\}\\
&&\qquad\quad\hspace*{11pt}{}\times f(y,\Sigma).
\end{eqnarray*}
\end{lemma}

\begin{lemma}\label{matrix-norm}
Let $\{\Sigma_T\}\subset\mathbb{R}^{(K+1)\times(K+1)}$ be a
sequence of positive definite matrices with $K+1\leq T$. If $K$ is
fixed with respect to $T$ and $\|\Sigma_T-\Sigma\|_{2}=O(1/T)$ for a
positive definite matrix $\Sigma$, then we have
\[
\bigl\|\Sigma_T^{-1}-\Sigma^{-1}\bigr\|_2=O(1/T).
\]
\end{lemma}

\begin{lemma}\label{remainder}
Let $\tilde{\Sigma}_T(y)$ be a $(K+1)\times(K+1)$ positive symmetric
matrix which depends on $y\in\mathbb{R}^{K+1}$. Assume that
$\sup_{y\in\mathbb{R}^{K+1}}\|\tilde{\Sigma}_T(y)-\Sigma
\|_{2}\leq
\|\Sigma_T-\Sigma\|_{2}=O(1/T)$ for a positive definite matrix
$\Sigma$. Let $R_T=\Sigma_T-\Sigma$. If $K$ is fixed with respect to
$T$, we have
\[
\int_{y\in\mathbb{R}^{K+1}}\biggl\llvert \operatorname{vec}(R_T)'
\frac{\partial
^2 f}{\partial
\operatorname{vec}(\Sigma)\operatorname{vec}(\Sigma)}\bigl(y,\tilde{\Sigma}_T(y)\bigr)
\operatorname{vec}(R_T)\biggr\rrvert \,dy=O\bigl(1/T^2
\bigr).
\]
\end{lemma}

\begin{pf*}{Proof of Theorem~\ref{t-stat}}
For the convenience of our presentation, we ignore the functional
symmetry of the covariance matrix in the proof. With some proper
modifications, we can extend the results to the case where the
functional symmetry is taken into consideration. Let
$|\mathcal{M}_1|=|\mathcal{M}_2|=\cdots=|\mathcal{M}_K|=q$. Define
$Y_i=\sqrt{q}(\hat{\mu}_i-\mu_0)$, and
$\bar{Y}=\frac{1}{K}\sum^{K}_{i=1}Y_i$ and
$S^2_Y=\frac{1}{K-1}\sum^{K}_{i=1}(Y_i-\bar{Y})^2$ as the sample
mean and sample variance of $\{Y_i\}^{K}_{i=1}$, respectively. Note
that $T_K(Y)=\sqrt{K}\bar{Y}/S_Y$, where $Y=(Y_1,Y_2,\ldots,Y_K)'$.
Simple algebra yields that
\[
\sigma_{ij}:=\operatorname{Cov}(Y_i,Y_j)=\sum
^{q-1}_{h=1-q} \biggl(\frac
{q-|h|}{q} \biggr)
\gamma_X\bigl(h-(j-i)q\bigr).
\]
Notice that $Y$ follows a normal distribution with mean zero and
covariance matrix~$\Sigma_T$, where
$\Sigma_T=(\sigma_{ij})_{i,j=1}^K$. The density function of $Y$ is
given by
\[
f(y,\Sigma_T)=(2\pi)^{-K/2}|\Sigma_T|^{-1/2}
\exp \bigl(-\tfrac
{1}{2}y'\Sigma_T^{-1}y
\bigr).
\]
Under the assumption
$\sum^{+\infty}_{h=-\infty}h^2|\gamma_X(h)|<\infty$, it is
straightforward to see that $\|\Sigma_T-\sigma^2I_K\|_2=O(1/T)$.
Taking a Taylor expansion of $f(y,\Sigma_T)$ around elements of the
matrix $\sigma^2I_K$, we have
\begin{eqnarray*}
f(y,\Sigma_T)&=&f\bigl(y,\sigma^2I_K\bigr)+
\biggl\{\frac{\partial f}{\partial
\operatorname{vec}(\Sigma)}\bigl(y,\sigma^2I_K\bigr)
\biggr\}'\operatorname{vec}\bigl(\Sigma _T-
\sigma^2I_K\bigr)
\\
&&{}+\operatorname{vec}\bigl(\Sigma_T-\sigma^2I_K
\bigr)'\frac{\partial^2 f}{\partial
\operatorname{vec}(\Sigma)\operatorname{vec}(\Sigma)}\bigl(y,\tilde{\Sigma}_T(y)
\bigr) \operatorname{vec}\bigl(\Sigma_T-\sigma^2I_K
\bigr),
\end{eqnarray*}
where
$\sup_{y\in\mathbb{R}^K}\|\tilde{\Sigma}_T(y)-\sigma
^2I_{K}\|_2\leq\|\Sigma_T-\sigma^2I_K\|_2=O(1/T)$.
By Lemmas~\ref{first-deri} and~\ref{remainder}, we get
\[
\frac{\partial f}{\partial\operatorname{vec}(\Sigma)} \bigl(y,\sigma^2I_K\bigr)=f\bigl(y,
\sigma^2I_K\bigr) \biggl\{-\frac{1}{2\sigma^2}
\operatorname{vec}(I_K)+\frac{1}{2\sigma^4}y\otimes y \biggr\}
\]
and
%
\begin{eqnarray}
\label{R-term}\qquad
&&\int_{y\in\mathbb{R}^K}\biggl\llvert \operatorname{vec}
\bigl(\Sigma_T-\sigma ^2I_K
\bigr)'\,\frac{\partial^2
f}{\partial
\operatorname{vec}(\Sigma)\operatorname{vec}(\Sigma)}\bigl(y,\tilde{\Sigma}_T(y)
\bigr) \operatorname{vec}\bigl(\Sigma_T-\sigma^2I_K
\bigr)\biggr\rrvert \,dy\nonumber\\[-8pt]\\[-8pt]
&&\qquad=O \biggl(\frac{1}{T^2} \biggr),\nonumber
\end{eqnarray}
which imply that
\begin{eqnarray*}
f(y,\Sigma_T)&=&f\bigl(y,\sigma^2I_K\bigr)
\Biggl\{1-\frac{1}{2\sigma^2}\sum^{K}_{i=1}
\bigl(\sigma_{ii}-\sigma^2\bigr) \Biggr\}
\\
&&{}+\frac{1}{2\sigma^4}f\bigl(y,\sigma^2I_K\bigr)\sum
^{K}_{i=1}\sum^{K}_{j=1}
\bigl(\sigma_{ij}-\sigma^2\delta_{ij}
\bigr)y_iy_j+R(y)
\\
&=&g\bigl(y,\sigma^2I_K\bigr)+R(y),
\end{eqnarray*}
where $g$ denotes the major term, $R(y)$ is the remainder term and
$\delta_{ij}=\mathbf{I}\{i=j\}$ is the Kronecker's delta. Define
$\tilde{\Psi}(x;K)=\int_{\{|T_K(y)|>x\}}g(y,\sigma^2I_K)\,dy$. By~(\ref{R-term}), we see that
\[
\sup_{x\in\mathbb{R}}\biggl\llvert \int_{\{|T_K(y)|>x\}}f(y,
\Sigma _T)\,dy-\tilde{\Psi}(x;K)\biggr\rrvert \leq \int
_{\mathbb{R}^K}\bigl|R(y)\bigr|\,dy=O\bigl(1/T^2\bigr).
\]
It follows from some simple calculation that
\[
\tilde{\Psi}(x;K)= \Biggl\{1-\frac{1}{2\sigma^2}\sum
^{K}_{i=1}\bigl(\sigma_{ii}-
\sigma^2\bigr) \Biggr\}P\bigl(|t_{K-1}|>x\bigr)+\frac
{1}{2\sigma^2}(J_1+J_2),
\]
where
\[
J_1=\sum^{K}_{i=1}\bigl(
\sigma_{ii}-\sigma^2\bigr)E\bigl[\mathbf{I}\bigl\{\bigl|\tilde
{T}_K(v)\bigr|>x\bigr\}v_i^2\bigr],\qquad\!\!
J_2=\sum_{i\neq j}\sigma_{ij}E
\bigl[\mathbf{I}\bigl\{ \bigl|\tilde{T}_K(v)\bigr|>x\bigr\}v_iv_j
\bigr].
\]
Here $\{v_i\}_{i=1}^{K}$ are i.i.d. standard normal random variables
and $\tilde{T}_K(v)=\break\sqrt{K}\bar{v}/S_v$ is the $t$ statistic based
on $\{v_i\}$ with $\bar{v}=\frac{1}{K}\sum^{K}_{i=1}v_i$ and
$S_v^2=\frac{1}{K-1}\sum^{K}_{i=1}(v_i-\bar{v})^2$. Let
$U=K\bar{v}^2$ and $D=(K-1)S^2_v$. Then $U\sim\chi^2_1$,
$D\sim\chi^2_{K-1}$ and $U$ and $D$ are independent. We define that
\begin{eqnarray*}
&&
E\bigl[\mathbf{I}\bigl\{\bigl|\tilde{T}_K(v)\bigr|>x\bigr
\}v_i^2\bigr]\\
&&\qquad=\frac{1}{K}E\Biggl[\mathbf{I}\bigl\{
\bigl|\tilde{T}_K(v)\bigr|>x\bigr\}\sum_{i=1}^{K}v_i^2
\Biggr]
\\
&&\qquad=\frac{1}{K}E\bigl[\mathbf{I}\bigl\{\bigl|\tilde{T}_K(v)\bigr|>x\bigr
\}U\bigr]+\frac
{1}{K}E\bigl[\mathbf{I}\bigl\{\bigl|\tilde{T}_K(v)\bigr|>x
\bigr\}D\bigr] 
\\
&&\qquad=\frac{1}{K}E \biggl[UG_{K-1} \biggl(\frac{(K-1)U}{x^2} \biggr)
\biggr]+\frac{1}{K}E \biggl[D-DG_{1} \biggl(\frac{Dx^2}{K-1}
\biggr) \biggr]
\end{eqnarray*}
and
\begin{eqnarray*}
&&
E\bigl[\mathbf{I}\bigl\{\bigl|\tilde{T}_K(v)\bigr|>x\bigr
\}v_iv_j\bigr]\\
&&\qquad=\frac{1}{K(K-1)}E\biggl[\mathbf {I}\bigl
\{\bigl|\tilde{T}_K(v)\bigr|>x\bigr\}\sum_{i\neq
j}v_iv_j
\biggr]
\\
&&\qquad=\frac{1}{K-1}E\bigl[\mathbf{I}\bigl\{\bigl|\tilde{T}_K(v)\bigr|>x\bigr
\}U\bigr]-\frac
{1}{K(K-1)}E\Biggl[\mathbf{I}\bigl\{\bigl|\tilde{T}_K(v)\bigr|>x
\bigr\}\sum_{i=1}^{K}v_i^2
\Biggr]
\\
&&\qquad=\frac{1}{K}E \biggl[UG_{K-1} \biggl(\frac{(K-1)U}{x^2} \biggr)
\biggr]-\frac{1}{K(K-1)}E \biggl[D-DG_{1} \biggl(\frac
{Dx^2}{K-1}
\biggr) \biggr].
\end{eqnarray*}
We then have
%
\begin{eqnarray}
\label{expan} P\bigl(|T_K|>x\bigr)&=&\tilde{\Psi}(x;K)+O\bigl(1/T^2
\bigr)
\nonumber
\\
&=& \{1-\alpha \}P\bigl(|t_{K-1}|>x\bigr)+\beta E \biggl[UG_{K-1} \biggl(
\frac{(K-1)U}{x^2} \biggr) \biggr]
\\
&&{}+\tau \biggl\{K-1-E \biggl[DG_{1} \biggl(\frac{Dx^2}{K-1} \biggr)
\biggr] \biggr\}+O\bigl(1/T^2\bigr),
\nonumber
\end{eqnarray}
uniformly for $x\in\mathbb{R}$, where the coefficients are given by
\begin{eqnarray*}
\alpha&=&\frac{1}{2\sigma^2}\sum^{K}_{i=1}
\bigl(\sigma_{ii}-\sigma ^2\bigr)=-\frac{K^2B}{2\sigma^2T}+O
\bigl(1/T^2\bigr),
\\
\beta&=&\frac{1}{2K\sigma^2}\sum^{K}_{i=1}\sum
^{K}_{j=1}\bigl(\sigma _{ij}-
\delta_{ij}\sigma^2\bigr)=-\frac{B}{2\sigma^2T}+O
\bigl(1/T^2\bigr)
\end{eqnarray*}
and
\[
\tau=\frac{1}{2K\sigma^2}\sum^{K}_{i=1}
\bigl(\sigma_{ii}-\sigma ^2\bigr)-\frac{1}{2K(K-1)\sigma^2}\sum
_{i\neq
j}\sigma_{ij}=-\frac{(K+1)B}{2\sigma^2T}+O
\bigl(1/T^2\bigr).
\]
The conclusion
thus follows from equation (\ref{expan}).
\end{pf*}

\section*{Acknowledgments}

The authors would like to thank the Associate Editor and the
reviewers for their constructive comments, which substantially
improve the paper.

\begin{supplement}
\stitle{Proofs of the other results in Sections
\ref{sechighorder}--\ref{secbootstrap} and simulation results.}
\slink[doi]{10.1214/13-AOS1113SUPP} 
\sdatatype{.pdf}
\sfilename{aos1113\_supp.pdf}
\sdescription{This supplement contains proofs of the other main
results in Sections~\ref{sechighorder}--\ref{secbootstrap} and
some simulation results.}
\end{supplement}


\printaddresses

\end{document}